\documentclass[11pt]{article}

\usepackage{amssymb,amsmath,amsfonts,amsthm}
\usepackage{latexsym}
\usepackage{graphics}
\usepackage{indentfirst}
\usepackage{hyperref}
\usepackage{comment}
\usepackage[shortlabels]{enumitem}
\usepackage[english]{babel}
\usepackage{tikz}
\usepackage{esdiff}

\usetikzlibrary{arrows.meta}

\newtheorem*{thm*}{Theorem}
\newtheorem{thm}{Theorem}

\newtheorem{conj}[thm]{Conjecture}

\newcommand{\N}{\mathbb{N}}

\allowdisplaybreaks

\begin{document}

\title{A Short Proof that the List Packing Number of any Graph is Well Defined}

\author{Jeffrey A. Mudrock$^1$}

\footnotetext[1]{Department of Mathematics, College of Lake County, Grayslake, IL 60030.  E-mail:  {\tt {jmudrock@clcillinois.edu}}}

\maketitle

\begin{abstract}

List packing is a notion that was introduced in 2021 (by Cambie et al.).  The list packing number of a graph $G$, denoted $\chi_{\ell}^*(G)$, is the least $k$ such that for any list assignment $L$ that assigns $k$ colors to each vertex of $G$, there is a set of $k$ proper $L$-colorings of $G$, $\{f_1, \ldots, f_k \}$, with the property $f_i(v) \neq f_j(v)$ whenever $1 \leq i < j \leq k$ and $v \in V(G)$.  We present a short proof that for any graph $G$, $\chi_{\ell}^*(G) \leq |V(G)|$.  Interestingly, our proof makes use of Galvin's celebrated result that the list chromatic number of the line graph of any bipartite multigraph equals its chromatic number.  

\medskip



\end{abstract}

\section{Introduction}\label{intro}

In this note all graphs are nonempty, finite, simple graphs unless otherwise noted.  Generally speaking we follow West~\cite{W01} for terminology and notation.  The set of natural numbers is $\N = \{1,2,3, \ldots \}$.  For $m \in \N$, we write $[m]$ for the set $\{1, \ldots, m \}$.  The chromatic number of a graph $G$ is denoted $\chi(G)$.   We write $K_{n,m}$ for complete bipartite graphs with partite sets of size $n$ and $m$.  The Cartesian product of graphs $G$ and $H$, denoted $G \square H$, is the graph with vertex set $V(G) \times V(H)$ and edges created so that $(u,v)$ is adjacent to $(u',v')$ if and only if either $u=u'$ and $vv' \in E(H)$ or $v=v'$ and $uu' \in E(G)$.  We will use the well-known fact that $\chi(G \square H) = \max \{\chi(G), \chi(H) \}$.  The line graph of $G$, denoted $L(G)$, is the graph with vertex set $E(G)$ such that distinct edges $e_1, e_2 \in E(G)$ are adjacent in $L(G)$ when $e_1$ and $e_2$ are incident in $G$.  Note that $L(K_{n,m}) = K_n \square K_m$.

\subsection{List Packing and the Main Theorem} \label{basic}

In classical vertex coloring we wish to color the vertices of a graph $G$ with up to $m$ colors from $[m]$ so that adjacent vertices receive different colors, a so-called \emph{proper $m$-coloring}.  List coloring is a well-known variation on classical vertex coloring that was introduced independently by Vizing~\cite{V76} and Erd\H{o}s, Rubin, and Taylor~\cite{ET79} in the 1970s.  For list coloring, we associate a \emph{list assignment} $L$ with a graph $G$ such that each vertex $v \in V(G)$ is assigned a list of colors $L(v)$ (we say $L$ is a list assignment for $G$).  Then, $G$ is \emph{$L$-colorable} if there exists a proper coloring $f$ of $G$ such that $f(v) \in L(v)$ for each $v \in V(G)$ (we refer to $f$ as a \emph{proper $L$-coloring} of $G$)~\footnote{A coloring $f$ of $G$ satisfying $f(v) \in L(v)$ for each $v \in V(G)$ that is not necessarily proper is called an \emph{$L$-coloring} of $G$}.  A list assignment $L$ for $G$ is said to be a \emph{$k$-assignment} if $|L(v)|=k$ for each $v \in V(G)$.  The \emph{list chromatic number} of $G$, denoted $\chi_{\ell}(G)$, is the smallest $k$ such that $G$ is $L$-colorable whenever $L$ is a $k$-assignment for $G$.  It is obvious that for any graph $G$, $\chi(G) \leq \chi_{\ell}(G)$.

One of the most famous open questions about list coloring is the List Coloring Conjecture which was formulated by many different researchers and has received considerable attention in the literature (see~\cite{HC92}).  
\begin{conj}[{\bf List Coloring Conjecture}] \label{conj: Edge List Coloring Conjecture}  If $G$ is a loopless multigraph, then $\chi(L(G)) = \chi_{\ell}(L(G))$. 
\end{conj}
In 1995, Galvin famously proved the following.
\begin{thm}[\cite{G95}] \label{thm: Galvin} If $G$ is a bipartite multigraph, then $\chi(L(G)) = \chi_{\ell}(L(G))$.
\end{thm}
Galvin's proof of Theorem~\ref{thm: Galvin} is regarded by many as one of the most beautiful proofs on the topic of list coloring (see~\cite{AZ18}).

List packing is a relatively new notion that was first suggested by Alon, Fellows, and Hare in 1996~\cite{AF96}.  This suggestion was not formally embraced until a recent paper of Cambie et al.~\cite{CC21}.  We now mention some important definitions.  Suppose $L$ is a list assignment for a graph $G$.  An \emph{$L$-packing of $G$ of size $k$} is a set of $k$ $L$-colorings of $G$, $\{f_1, \ldots, f_k \}$, such that $f_i(v) \neq f_j(v)$ whenever $i, j \in [k]$, $i \neq j$, and $v \in V(G)$.  Moreover, we say that $\{f_1, \ldots, f_k \}$ is \emph{proper} if $f_i$ is a proper $L$-coloring of $G$ for each $i \in [k]$. The \emph{list packing number} of $G$, denoted $\chi_{\ell}^*(G)$, is the least $k$ such that $G$ has a proper $L$-packing of size $k$ whenever $L$ is a $k$-assignment for $G$.  Clearly, for any graph $G$, $\chi(G) \leq \chi_{\ell}(G) \leq \chi_{\ell}^*(G)$.

In~\cite{CC21}, while defining $\chi_{\ell}^*(G)$, the authors remark that ``The reader might already find it interesting that such a minimal $k$ is well defined.".  The authors go on to show that for any graph $G$, $\chi_{\ell}^*(G) \leq |V(G)|$ and equality holds if and only if $G$ is complete, and they conjecture that there is a $C > 0$ such that $\chi_{\ell}^*(G) \leq C \chi_{\ell}(G)$. In this note we present a short proof that for any graph $G$, $\chi_{\ell}^*(G) \leq |V(G)|$ which makes use of Theorem~\ref{thm: Galvin}.  In particular, we prove the following.
\begin{thm} \label{thm: simpleproof}
$\chi_{\ell}^*(K_n) = n$.
\end{thm}
Note that Theorem~\ref{thm: simpleproof} implies that for any graph $G$, $\chi_{\ell}^*(G) \leq |V(G)|$ since $\chi_{\ell}^*(H) \leq \chi_{\ell}^*(G)$ whenever $H$ is a subgraph of $G$. 

\section{Proof of Theorem~\ref{thm: simpleproof}} \label{main}

\begin{proof}
Since $\chi(K_n) = n$, $\chi_{\ell}^*(K_n) \geq n$.  Suppose $G=K_n$ and $V(G) = \{v_1, \ldots, v_n\}$.  Let $L$ be an arbitrary $m$-assignment for $G$ with $m \geq n$.  To prove the desired result, we will show that there is a proper $L$-packing of $G$ of size $m$.

Suppose $H = G \square K_m$ and the vertices of the copy of $K_m$ used to form $H$ are $\{u_1, \ldots, u_m \}$.  Let $L_H$ be the $m$-assignment for $H$ given by $L_H(v_i,u_j) = L(v_i)$ for each~$(i,j) \in [n] \times [m]$.  By Theorem~\ref{thm: Galvin} and properties of the Cartesian product of graphs, $m = \chi(K_n \square K_m) = \chi(L(K_{n,m})) = \chi_{\ell}(L(K_{n,m})) = \chi_{\ell}(H)$.  Consequently, there is a proper $L_H$-coloring of $H$ which we will name $f$.  Now, for each $j \in [m]$ let $f_j$ be the proper $L$-coloring of $G$ given by $f_j(v_i) = f(v_i,u_j)$ for each $i \in [n]$.  Finally, notice that $\{f_1, \dots, f_m \}$ is a proper $L$-packing of $G$ of size $m$.
\end{proof}

{\bf Acknowledgment.}  The author would like to thank Hemanshu Kaul for helpful conversations regarding the contents of this note.  The author would also like to thank the anonymous referees for their helpful comments on this note.

\end{document}